\title{\vspace{-1cm} On the strong chromatic number of random graphs}

\author{
Po-Shen Loh \thanks{Department of Mathematics,
Princeton University, Princeton, NJ 08544. E-mail:
ploh@math.princeton.edu.
Research supported in part by a Fannie and John Hertz Foundation Fellowship, an NSF Graduate
Research Fellowship, and a Princeton Centennial Fellowship.}
\and Benny Sudakov \thanks{Department of Mathematics, Princeton University, Princeton, NJ 08544, and
Institute for Advanced Study, Princeton. E-mail:
bsudakov@math.princeton.edu.
Research supported in part by NSF CAREER award DMS-0546523, NSF grant
DMS-0355497, USA-Israeli BSF grant, Alfred P. Sloan fellowship, and
the State of New Jersey.
}
}

\documentclass[11pt]{article}
\usepackage{amsmath, amssymb, epsfig, graphicx, enumerate}

\newtheorem{theorem}{Theorem}[section]

\newtheorem{lemma}[theorem]{Lemma}
\newtheorem{proposition}[theorem]{Proposition}

\newcommand{\ep}{\epsilon}
\newcommand{\pr}[1]{\mathbb{P}\left[#1\right]}

\newcommand{\gnp}{G_{n, p}}
\newcommand{\numparts}{{\big\lceil\frac{1}{p}\big\rceil}}
\newcommand{\smallnumparts}{{\lceil 1/p \rceil}}
\newcommand{\twicenumparts}{{2\lceil\frac{1}{p}\rceil}}

\date{}
\begin{document}
\maketitle

\begin{abstract}
  Let $G$ be a graph with $n$ vertices, and let $k$ be an integer dividing $n$.  $G$ is
  said to be strongly $k$-colorable if for every partition of $V(G)$ into disjoint sets
  $V_1 \cup \ldots \cup V_r$, all of size exactly $k$, there exists a proper vertex
  $k$-coloring of $G$ with each color appearing exactly once in each $V_i$.  In the case
  when $k$ does not divide $n$, $G$ is defined to be strongly $k$-colorable if the graph
  obtained by adding $k \big\lceil \frac{n}{k} \big\rceil - n$ isolated vertices is
  strongly $k$-colorable.  The strong chromatic number of $G$ is the minimum $k$ for which
  $G$ is strongly $k$-colorable.  In this paper, we study the behavior of this parameter for
  the random graph $\gnp$.  In the dense case when $p\gg n^{-1/3}$, we prove that the strong 
  chromatic number is a.s.\ concentrated on one value $\Delta+1$, where $\Delta$ is the maximum 
  degree of the graph. We also obtain several weaker results for sparse random graphs.
\end{abstract}

\section{Introduction}

Let $G$ be a graph, and let $V_1$, \ldots, $V_r$ be disjoint subsets of its vertex set.
An \emph{independent transversal}\/ with respect to $\{V_i\}_{i=1}^r$ is an independent
set in $G$ which contains exactly one vertex from each $V_i$. The problem of finding
sufficient conditions for the existence of an independent transversal, in terms of the
ratio between the part sizes and the maximum degree $\Delta$ of the graph, dates back to
1975, when it was raised by Bollob\'as, Erd\H{o}s, and Szemer\'edi \cite{BES}.  Since
then, much work has been done \cite{ABZ, Al-problems-results, Al-strong-chromatic,
  Ha-note, Ha-strong-chromatic, HS, Ji, Me, ST, Yu}, and this basic concept has also
appeared in several other contexts, such as linear arboricity \cite{Al-linear-arboricity},
vertex list coloring \cite{Re, RS, BH}, and cooperative coloring \cite{AH, LS}.  In the
general case, it was proved by Haxell \cite{Ha-note} that an independent transversal
exists as long as all parts have size at least $2\Delta$. The sharpness of this bound was
shown by Szab\'o and Tardos \cite{ST}, extending earlier results of \cite{Ji} and \cite{Yu}.
On the other hand, we proved in \cite{LS} that the upper bound can be further reduced
to $(1+o(1))\Delta$ if no vertex has more than $o(\Delta)$ neighbors in any single part.
Such a condition arises naturally in certain applications, e.g., vertex list coloring.

In the case when all of the $V_i$ are of the same size $k$, it is natural to ask when it
is possible to find not just one, but $k$ disjoint independent transversals with respect
to the $\{V_i\}$.  This is closely related to the following notion of strong colorability. Given a
graph $G$ with $n$ vertices and a positive integer $k$ dividing $n$, we say that $G$ is
\emph{strongly $k$-colorable}\/ if for every partition of $V(G)$ into disjoint sets $V_1
\cup \ldots \cup V_r$, all of size exactly $k$, there exists a proper vertex $k$-coloring of
$G$ with each color appearing exactly once in each $V_i$. Notice that $G$ is strongly 
$k$-colorable iff the chromatic number of any graph obtained from $G$ by adding a union of 
vertex disjoint $k$-cliques is $k$. If $k$ does not divide $n$,
then we say that $G$ is strongly $k$-colorable if the graph obtained by adding $k
\big\lceil \frac{n}{k} \big\rceil - n$ isolated vertices is strongly $k$-colorable.  The
\emph{strong chromatic number}\/ of $G$, denoted $s\chi(G)$, is the minimum $k$ for which
$G$ is strongly $k$-colorable.

The concept of strong chromatic number first appeared independently in work by Alon
\cite{Al-linear-arboricity} and Fellows \cite{Fe}.  It was also the crux of the
longstanding ``cycle plus triangles'' problem popularized by Erd\H{o}s, which was to
show that the strong chromatic number of the cycle on $3n$ vertices is three. That
problem was solved by Fleischner and Stiebitz \cite{FS-c3n}.  The strong
chromatic number is known \cite{Fe} to be monotonic in the sense that strong
$k$-colorability implies strong $(k+1)$-colorability. It is also easy to see that $s\chi(G)$
must always be strictly greater than the maximum degree $\Delta$:
simply take $V_1$ to be the neighborhood of a vertex of maximal degree, 
and partition the rest of the vertices arbitrarily. The intriguing question of bounding the strong chromatic 
number in terms of the maximum degree has not yet been answered
completely. Alon \cite{Al-strong-chromatic} showed that there exists a constant $c$ such
that $s\chi \leq c\Delta$ for every graph.  Later, Haxell \cite{Ha-strong-chromatic}
improved the bound by showing that it is enough to use $c=3$, and in fact even $c=3-\ep$
for $\ep$ up to $1/4$ \cite{Ha-private-comm}.  On the other hand, Fleischner and Stiebitz \cite{FS-2d}
observed that the disjoint union of complete bipartite graphs $K_{\Delta, \Delta}$ cannot be strongly $(2\Delta-1)$-colored.
Indeed, put each part of one of the $K_{\Delta, \Delta}$ into the sets $V_1$ and $V_2$, respectively.
Then these $2\Delta$ vertices should get different colors. It is believed that this lower bound is tight and 
the strong chromatic number of any graph with maximum degree $\Delta$ should be at most $2\Delta$.

It is natural to wonder what is the asymptotic behavior of the strong chromatic number for
the random graph $\gnp$, relative to the maximum degree of the graph. As usual, $\gnp$ is
the probability space of all labeled graphs on $n$ vertices, where every edge appears
randomly and independently with probability $p=p(n)$. We say that the random graph
possesses a graph property $\cal P$ {\em almost surely}, or a.s.\ for brevity, if the
probability that $\gnp$ satisfies $\cal P$ tends to 1 as $n$ tends to infinity. One of the
most interesting phenomena discovered in the study of random graphs is that many natural
graph invariants are highly concentrated (see, e.g., \cite {Ma} for the result on the
clique number and \cite{SS, Lu, AK} for the concentration of the chromatic number). In
this paper we show that the strong chromatic number is another example of a tightly
concentrated graph parameter.  For dense random graphs, it turns out that we can
concentrate $s\chi(\gnp)$ on a single value, and for some smaller values of $p$ we were
only able to determine $s\chi(\gnp)$ asymptotically. In the statement of our first result,
and in the rest of this paper, the notation $f(n) \gg g(n)$ means that $f/g \rightarrow
\infty$ together with $n$.  Also, all logarithms are in the natural base $e$.

\begin{theorem}
  \label{thm:main}
  Let $\Delta$ be the maximum degree of the random graph $\gnp$, where $p < 1-\theta$
  for any arbitrary constant $\theta > 0$.
  \begin{description}
  \item[\hspace{12pt}\textbf{(i)}] If $p \gg \left(\frac{\log^4 n}{n}\right)^{1/3}$, then
almost surely the strong chromatic
    number of $\gnp$ equals $\Delta+1$.
\item[\hspace{12pt}\textbf{(ii)}] If $p \gg \left(\frac{\log n}{n}\right)^{1/2}$, then a.s.\ the strong chromatic
    number of $\gnp$ is $(1+o(1))\Delta$.
  \end{description}
\end{theorem}

Unfortunately, our approach breaks down completely when $p \ll n^{-1/2}$. However, for this range 
of $p$, we have a different argument which shows how to find  at least one independent transversal.

\begin{theorem}
\label{thm:main2}
Let $\Delta$ be the maximum degree of the random graph $\gnp$. If
$p \geq \frac{\log^4 n}{n}$,
then almost surely every collection of disjoint subsets $V_1, \ldots, V_r$ of $\gnp$ with all
$|V_i|\geq (1+o(1))\Delta$ has an independent transversal.
\end{theorem}

This rest of this paper is organized as follows.  In Section \ref{sec:first-two}, we prove
both parts of our first theorem concerning the strong chromatic number of relatively dense
random graphs.  We then shift our attention to the sparser case, proving our second result
about transversals in Section \ref{sec:third}.  The last section of our paper contains
some concluding remarks.  Throughout this exposition, we will make no attempt to optimize
absolute constants, and will often omit floor and ceiling signs whenever they are not
crucial, for the sake of clarity of presentation.

\section{Strong chromatic number}
\label{sec:first-two}

In this section, we prove Theorem \ref{thm:main}, which determines the value of the strong
chromatic number of a rather dense random graph. To this end, we first prove several lemmas
that establish certain useful properties of random graphs. We will use these properties
to find a partition of $\gnp$ into independent transversals.

\subsection{Properties of random graphs}
\label{sec:whp}

\begin{lemma}
  \label{lem:degree-sequence}
  Let $\theta > 0$ be an arbitrary fixed constant.
  If $\sqrt{\frac{\log n}{n}} \ll p < 1-\theta$ then a.s.\ $\gnp$ has the
  following properties.
  \begin{description}
  \item[\hspace{12pt}\textbf{(i)}] No pair of distinct vertices has more than $(1+o(1))np^2$ common neighbors.
  \item[\hspace{12pt}\textbf{(ii)}] The maximum degree is strictly between $np$ and $1.01np$, and there is a unique
    vertex of maximum degree.
  \item[\hspace{12pt}\textbf{(iii)}] The gap between the maximum degree and the next largest degree is at least
    $\frac{\sqrt{np}}{\log n}$.
  \end{description}
\end{lemma}

\noindent {\bf Proof.}\, For the first property, fix an arbitrary constant $\delta>0$ and two 
distinct vertices $u$ and $v$.
Their codegree $X$ is binomially distributed with parameters $n-2$ and $p^2$.  Thus by the
Chernoff bound (see, e.g., Appendix A in \cite{AS}), $\pr{X \geq (1+\delta)np^2} \leq 
e^{-\Theta(\delta^2np^2)}=
o(n^{-2})$.  Taking a union bound over all $O(n^2)$ choices for $u$ and $v$, we find that
the probability that the first property is not satisfied tends to 0 as $n \rightarrow
\infty$.  The second and third claims are special cases of Corollary 3.13 and Theorem 3.15
in \cite{Bol}, respectively.  \hfill $\Box$

\begin{lemma}
  \label{lem:dominate}
  Let $\alpha > 0$ be an arbitrary fixed constant and let $\sqrt{\frac{\log n}{n}} \ll p
  \leq \frac{3}{5}$. Then almost surely $\gnp$ does not contain a set $U$ of size $\alpha np$
  and $50 \log n$ sets $T_i$, $|T_i| \leq \numparts$, such that all the sets are disjoint and
  for every $i$ all but at most $\alpha np/50$ vertices in $U$ have neighbors in $T_i$.
\end{lemma}

\noindent {\bf Proof.}\, Fix sets $U$ and $\{T_i\}$ as specified above. If all but at most
$\alpha np/50$ vertices in $U$ have neighbors in $T_i$, we say for brevity that $T_i$ {\em
  almost dominates $U$}. For a given vertex $v$, the probability that it has a neighbor in
$T_i$ is $1 - (1-p)^{|T_i|}\leq 1 - (1-p)^{\smallnumparts}<7/8$ for all $p \leq 3/5$,
since $1 - (1-p)^{\smallnumparts}$ is maximal in that range when $p \rightarrow 1/2$ from below.
Therefore, by a union bound we have

\begin{eqnarray*}
\pr{\text{$T_i$ almost dominates $U$}} 
  &\leq& {\alpha n p \choose {\alpha np-\alpha n p/50}} \left(\frac{7}{8}\right)^{\alpha np-\alpha n 
p/50} \ = \ 
{\alpha n p \choose \alpha n p/50}\left(\frac{7}{8}\right)^{49 \alpha np/50}\\
&\leq& \left( 50e \Big(\frac{7}{8}\Big)^{49}\right)^{\alpha np/50} \ < \ 3^{-\alpha np/50}.
\end{eqnarray*}
Since all sets $T_i$ are disjoint, the events that $T_i$ and $T_j$
almost dominate $U$ are independent. 
This implies that
\begin{displaymath}
  \pr{\text{every $T_i$ almost dominates $U$}} \ \leq \ \left(3^{-\alpha np/50}\right)^{50
    \log n} \ = \ 3^{-\alpha np\log n}.
\end{displaymath}
Using that $\log n/p =o(np)$ and $\smallnumparts\leq 2/p$, we can bound the probability
that there is a choice of $\{T_i\}$ and $U$ which violates the assertion of the lemma by
\begin{eqnarray*}
  \mathbb{P} &\leq& {n \choose \alpha np}
  \left[\frac{2}{p} {n \choose 2/p}\right]^{50 \log n} 3^{-\alpha np\log n}\\
  &\leq& n^{\alpha np} \left(\frac{2}{p}\right)^{50 \log n} n^{\frac{100 \log n}{p}}
3^{-\alpha np\log n}\\
  &=& e^{(1+o(1))\alpha np \log n} \cdot 3^{-\alpha np\log n} \ = \ o(1),
\end{eqnarray*}
so we are done.  \hfill
$\Box$

\begin{lemma}
  \label{lem:indep-trans}
Let $\alpha > 0$ be an arbitrary fixed constant and let
$\sqrt{\frac{\log n}{n}} \ll p \leq \frac{3}{5}$. Then almost surely every collection of at most $\numparts$
disjoint subsets of size $\alpha np$ in $\gnp$ has an independent transversal.
\end{lemma}

\noindent {\bf Proof.}\, Fix a collection of disjoint subsets $V_1, \dots, V_r$, $r \leq
\numparts$, of $\gnp$, each of size $\alpha np$.  A partial independent transversal $T$ is
an independent set with at most one vertex in every $V_i$, and we say that it almost
dominates some part if all but at most $\alpha np/50$ vertices in that part have neighbors
in $T$.  For every $V_i$, let $\{T_{ij}\}$ be a maximal collection of pairwise disjoint
partial independent transversals, each of which almost dominates $V_i$. Then, by Lemma
\ref{lem:dominate}, a.s.\ the total number of $T_{ij}$ must be at most $r (50 \log n)$.
Delete all the sets $T_{ij}$ from the graph, and let $\{V_i'\}$ be the remaining parts.
Clearly, it suffices to find an independent transversal among the $\{V_i'\}$.

Since $\log n/p =o(np)$ and each $T_{ij}$ is a partial transversal, each part loses a
total of $\leq r(50 \log n) \leq 50 \numparts \log n =o(np)$ vertices from the deletions.  We
can now use the greedy algorithm to find an independent transversal.  Take $v_1$ to be any
remaining vertex in $V'_1$, and iterate as follows.  Suppose that we already have
constructed a partial independent transversal $\{v_1, \ldots, v_{\ell-1}\}$ such that $v_i
\in V'_i$ for all $i < \ell$.  This partial independent transversal does not almost
dominate $V_\ell$, or else it would contradict the maximality of $\{T_{\ell j}\}$ above.
So, there are at least $\alpha np/50$ choices for $v_\ell \in V_\ell$ that would extend
the partial independent transversal $\{v_1, \ldots, v_{\ell-1}\}$.  Yet $V_\ell$ lost only
$o(np)$ vertices in the deletion process, so there is still a positive number of choices
for $v_\ell \in V_\ell'$ as well.  Proceeding in this way, we find a complete
independent transversal. \hfill $\Box$

\begin{lemma}
\label{lem:hall}
Let $\sqrt{\frac{\log n}{n}} \ll p \leq \frac{3}{5}$.  Then the following statement holds almost
surely.  For every choice of $s$ and $t$ that satisfies $np/2 \leq s \leq 2np$ and $40
\log n \leq t \leq s - 40 \numparts \log n$, $\gnp$ does not contain a collection of
disjoint subsets $U, T_1, \ldots, T_t$ such that $|U| = s$, each of the $|T_i| \leq
\numparts$, and at least $s-t$ vertices of $U$ have neighbors in every $T_i$.
\end{lemma}

\noindent {\bf Proof.}\, Fix some $(s, t)$ within the above range. As we saw in the proof of Lemma
\ref{lem:dominate}, for a given vertex $v$ the probability that it has a neighbor in $T_i$
is $1 - (1-p)^{|T_i|}\leq 1 - (1-p)^{\smallnumparts}<7/8$, and by disjointness 
these events are independent for all $1 \leq i \leq t$. Therefore
we can bound the the probability that there is a collection of sets which satisfies the 
above condition by
\begin{eqnarray}
\label{eq1}
\mathbb{P} &\leq&  {n \choose s} \left[\frac{2}{p} {n \choose 2/p}\right]^t 2^s 
\left(\frac{7}{8}\right)^{(s - t)t} \nonumber\\
&\leq&\frac{n^s}{s!} \, \big(n^{2/p}\big)^t\, 2^s \left(\frac{7}{8}\right)^{(s - t)t} \nonumber\\
&\leq& n^{s+2t/p} \left(\frac{7}{8}\right)^{(s - t)t} \,.
\end{eqnarray}
Throughout this bound, we use $\numparts \leq \frac{2}{p}$.  The first binomial
coefficient and the quantity in the square brackets bound the number of ways to choose the
sets $U$ and $\{T_i\}$. The $2^s$ bounds the number of ways to select a subset of size
$s -t$ from $U$, and the final factor bounds
the probability that all vertices in this subset have neighbors in every $T_i$.

The logarithm of (\ref{eq1}) is quadratic in $t$ with positive $t^2$-coefficient.
Therefore, the right hand side of (\ref{eq1}) is largest when $t$ is minimum or maximum in
its range $40 \log n \leq t \leq s-40\numparts \log n$.  Let us begin with the
small end, i.e., $t = 40 \log n$.  Then, since $\log n/p \ll np$ and $s\geq np/2$, we have
that
\begin{eqnarray*}
 n^{s+2t/p} \left(\frac{7}{8}\right)^{(s - t)t} 
&\leq& e^{(1+o(1))s\log n} \, \left(\frac{7}{8}\right)^{(40-o(1))s \log n}\\
&\leq&
e^{(1+o(1))s\log n} \,\, e^{-(4-o(1))s\log n} \ = \ o\big(n^{-2}\big).
\end{eqnarray*}
Similarly, if $t = s - 40\numparts \log n$, the bound is
\begin{eqnarray*}
n^{s+2t/p} \left(\frac{7}{8}\right)^{(s - t)t}
&\leq&  e^{3s \log n/p} \left(\frac{7}{8}\right)^
{(40-o(1))s\numparts \log n}\\ & \leq&
e^{3s \log n/p} \,\, e^{-(4-o(1))s\numparts \log n} \ = \ o\big(n^{-2}\big).
\end{eqnarray*}
Since the number of choices for $t$ and $s$ is at most $n^2$, we conclude that the
probability that the assertion of the lemma is violated is $o(1)$.
\hfill $\Box$

\subsection{Proof of Theorem \ref{thm:main}}
\label{sec:condition}

We start by proving part (i) of Theorem \ref{thm:main}.  If $\Delta$ is the maximum degree
of $\gnp$, then the strong chromatic number must be at least $\Delta+1$, as we already
mentioned in the introduction. Suppose that $G$ is a graph obtained from $\gnp$ by adding
$(\Delta+1)\lceil \frac{n}{\Delta+1}\rceil-n$ isolated vertices, and we have a partition
of $V(G)$ into $V_1 \cup \ldots \cup V_r$ with every $|V_i| = \Delta+1$. By Lemma
\ref{lem:degree-sequence}, $\Delta \geq np$ almost surely, so this implies that $r \leq
\numparts$.  Note that if $3/5 \leq p < 1-\theta$, then $r \leq 2$ and the theorem is an
immediate consequence of the following lemma.

\begin{lemma}
  \label{lem:p>3/5}
  Let $3/5 \leq p < 1-\theta$, where $\theta > 0$ is an arbitrary fixed constant, and let
  $V(G)=V_1 \cup V_2$ be a partition of the vertices of $G$ described above, with $|V_1| =
  |V_2| = \Delta+1$. Then a.s.\ $V_1$ can be perfectly matched to $V_2$ via non-edges of
  $G$.
\end{lemma}

\noindent {\bf Proof.}\, Without loss of generality, we may assume that $V_1$ contains at
most $n/2$ original vertices of $\gnp$.  Let $B \subset V_1$ be those original vertices.
The rest of $V_1$ consists of isolated vertices, so any perfect matching of $B$ to $V_2$
trivially extends to a full perfect matching between $V_1$ and $V_2$.  Therefore, by Hall's
theorem, it suffices to verify that each subset $A \subset B$ has at least $|A|$
non-neighbors in $V_2$.  If $A = \{v\}$ is a single vertex, this is immediate because
$|V_2| > \Delta \geq d(v)$.  For larger $A$, the Hall condition translates into
checking that $\Delta+1 - |N(A)| \geq |A|$, where $N(A)$ denotes the set of common
neighbors of $A$ in $V_2$.  Since $|A| \geq 2$ we have, by Lemma
\ref{lem:degree-sequence}(i), that the size of $N(A)$ is at most $(1+o(1))np^2$. So the
Hall condition is satisfied for all $A$ with $2 \leq |A| \leq \theta np/2 < \Delta -
(1+o(1))np^2$.

Let $c$ be a constant for which $p- 2p^c > 1/2$ for all $p$ in the range $[3/5, 1-\theta)$.  
One can easily show using a Chernoff bound
that a.s.\ every set of $c$ distinct vertices in $\gnp$ has at most
$2np^c$ common neighbors.  This implies that the Hall condition is also satisfied for all
$A$ of size at least $c$, since then
\begin{displaymath}
  \Delta+1 - |N(A)| \ > \ np - 2np^c \ > \ n/2 \ \geq \ |B| \ \geq \ |A|.
\end{displaymath}
Together with the previous paragraph, this completes the proof.
\hfill $\Box$

\vspace{0.35cm}

It remains to consider $p < 3/5$, so we will assume that bound on $p$ for the 
remainder of this section.  We use the following strategy to produce a partition of $\cup V_i$ 
into a disjoint union of independent transversals.

\begin{enumerate}
\item Find an independent transversal through the unique vertex of maximum degree $\Delta$,
  and delete this transversal from the graph.

\item As long as there exists a vertex $v$ which has at least $0.9 np$ neighbors in some part 
$V_i$, find an independent transversal $T$ through $v$, and delete $T$ from the graph.

\item As long as there exists a minimal partial independent transversal $T$ such that
all but at most $np/100$ vertices in some part $V_i$ have neighbors in $T$,
split $T$ into two nonempty ($|T| \geq 2$ because of Step 2) disjoint partial independent 
transversals $T_1 \cup T_2$.  Note that by minimality of $T$, each part $V_i$ contains a 
subset $U_i$ of at 
least $np/100$ vertices which have no neighbors in $T_1$. By Lemma \ref{lem:indep-trans}, there is 
an independent transversal through $\{U_i\}$, which can be used to extend $T_1$ to a full 
independent transversal $T'_1$. Delete $T'_1$ from the graph, and then perform the same
completion/deletion procedure for $T_2$.

\item Finally, we construct the rest of the independent transversals, building them 
simultaneously from $V_1$ to $V_r$ using Hall's
matching theorem.  Our deletions in Steps 1--3, together with the properties of $\gnp$ 
which we established in the previous subsection, will ensure that this is possible.
\end{enumerate}

The following lemma, which we prove later, ensures that we will
indeed find the independent transversals claimed in Steps 1--2.

\begin{lemma}
\label{lem:delete-single-vertices}
Let $V_1 \cup \ldots \cup V_r$ be the above partition of $V(G)$, and let $x$ be any vertex in this
graph.  
\begin{itemize}
\item
If $x$ is the unique vertex of maximum degree $\Delta$, then $G$ contains an independent
transversal through $x$. 
\item
If $x$ is not of maximum degree, then for all $k \leq \numparts$ and for any collection of subsets 
$V'_i \subset V_i$, $|V'_i|=\Delta+1 - k$, one of which contains $x$,
there exists an independent transversal through $x$ with respect to $\{V'_i\}$.
\end{itemize}
\end{lemma}

Let us bound the number of independent transversals we delete in the first 3 steps.
Note that if two vertices have at least $0.9np$ neighbors in the same $V_i$, since by
Lemma \ref{lem:degree-sequence} $|V_i| \leq \Delta+1 \leq 1.01np$, their codegree will be
at least $0.79 np \geq 1.01 np^2$, contradicting Lemma \ref{lem:degree-sequence}.
Therefore, during the first two steps, we will delete at most $r+1 \leq \numparts+1$
transversals.  Next, suppose that after deleting $O\big(\numparts \log n\big)$
independent transversals from $G$, we have that for some set $T$ all but at most $np/100$
vertices of some $V_i$ have neighbors in $T$.  Since $\numparts \log n \ll np$, this certainly
implies that the number of vertices in the original $V_i$ with no neighbors in $T$ was
bounded by $np/50$. Together with Lemma \ref{lem:dominate}, this ensures that for each fixed
$V_i$, $1 \leq i \leq r$, we never repeat Step 3 more than $50 \log n$ times.  Since each
iteration deletes two independent transversals and $r \leq
\numparts$, we conclude that by the time we reach Step 4, we have deleted at most
$1+\numparts + 100\numparts \log n <110\numparts \log n$ independent transversals from
$G$.

Let us now describe Step 4 in more detail.  At this point, all parts $V_i$ have the same size 
$|V_i| = s=\Delta+1-k$, where $k < 110\numparts \log n=o(np)$ is the total number of independent 
transversals 
deleted so far.  We build the remaining $s$ disjoint independent transversals simultaneously as 
follows. Start $s$ partial independent transversals $\{T_i\}_{i=1}^s$ by arbitrarily putting one vertex 
of $V_1$ 
into each $T_i$. Now suppose we already have disjoint partial independent transversals $\{T_i\}_{i=1}^s$
through $V_1, \ldots, V_\ell$. 
Create an auxiliary bipartite graph $H$ whose right side is $V_{\ell+1}$
and left side has $s$ vertices, identified with the transversals 
$\{T_i\}$.  Join the $i$-th vertex on the left side 
with a vertex $v\in V_{\ell+1}$ if and only if 
$v$ has no neighbors in $T_i$.  Then, a perfect matching in this graph will yield a simultaneous 
extension 
of each $T_i$ which covers $V_{\ell+1}$.

We ensure a perfect matching in $H$ by verifying the Hall condition, i.e.,  we show that for every
$t \leq s$, every set of $t$ vertices on the left side of $H$ has neighborhood on
the right side of size at least $t$.  Observe that after Step 3, for
every $T_i$ there are more than $np/100$ vertices in $V_{\ell+1}$ which have no neighbors in $T_i$.
Therefore every vertex on the left
side of $H$ has degree greater than $np/100$ and hence the Hall condition is
trivially satisfied for all $t \leq np/100$. If 
the Hall condition fails for some $np/100 <t \leq s-40\numparts \log n$, then 
by definition of $H$, there are $t$ partial independent transversals
among $\{T_i\}$ and a subset $W$ of $V_{\ell+1}$ of size  greater than $s-t$ such that every vertex of $W$ 
has 
neighbors 
in every one of these transversals (i.e., is not adjacent to them in $H$).
This contradicts Lemma \ref{lem:hall}, so the Hall condition also holds for these $t$.
It remains to check the case when $t > s - 40 \numparts \log n$.  Note
that given any vertex $v$ in $V_{\ell+1}$ and any collection of disjoint partial independent 
transversals,
the number of them in which $v$ can have a neighbor is at most the degree of $v$.  However,
we deleted the maximum degree vertex in Step 1, so by Lemma \ref{lem:degree-sequence} 
$d(v) \leq \Delta - \frac{\sqrt{np}}{\log n}$. Since
$p \gg \left(\frac{\log^4n}{n}\right)^{1/3}$,  this is less than 
$\Delta -150\numparts \log n \leq s - 40 \numparts \log n$.
Therefore, in the auxiliary graph $H$, any set of $t > s -40\numparts \log n$ vertices on the left side 
has 
neighborhood equal to the entire right side.  
Hence Hall's condition is satisfied for all $t$ and we can extend our transversals.
This completes the proof, since one can iterate this extension procedure to convert
all $T_i$ into full independent transversals. \hfill $\Box$

\vspace{0.35cm}

\noindent {\bf Proof of Lemma \ref{lem:delete-single-vertices}.}\, 
First, consider the case when $x$ is not the vertex of maximum degree $\Delta$ and we have
a collection of subsets $V'_i \subset V_i$ of size $\Delta+1-k$, where $k \leq \numparts$.
Without loss of generality, assume that $x \in V'_1$, and recall that by Lemma
\ref{lem:degree-sequence}, the maximum degree $\Delta$ satisfies $np < \Delta < 1.01np$.
If the number of neighbors of $x$ in every set $V'_i$, $i \geq 2$, is at most $0.96np$
then delete them and denote the resulting sets $V''_i$.  Since each $V''_i$ still has size
at least $\Delta+1-\numparts-0.96np>0.03np$, by Lemma \ref{lem:indep-trans} there exists a
partial independent transversal through $V''_2, \ldots, V''_r$, which together with $x$
provides a full independent transversal containing $x$.  Next, suppose that $x$ has at
least $0.96np$ neighbors in some part, say $V'_2$.  Since the degree of $x$ is less than
$\Delta < 1.01np$, it must then have less than $0.05np$ neighbors in every other $V'_i$.
Furthermore, since $x$ is not of maximum degree and $p \gg \left(\frac{\log^4
    n}{n}\right)^{1/3}$, Lemma \ref{lem:degree-sequence} implies that $(\Delta+1) -
d(x) \geq \frac{\sqrt{np}}{\log n} \gg \twicenumparts \geq r + k$.  Therefore
there are more than $r$ vertices in $V'_2$ not adjacent to $x$.  Also by Lemma
\ref{lem:degree-sequence}, the codegree of every pair of vertices is at most
$1.01np^2<0.61np$, so in particular no two vertices can both have $\geq 0.9np$ neighbors
in any given $V'_i$.  By the pigeonhole principle, there must be a vertex $y \in V'_2$ not
adjacent to $x$ with less than $0.9np$ neighbors in each of the other $V'_i$.  That means
that every other part has less than $0.05np$ neighbors of $x$ and $0.9np$ neighbors of $y$.
Since $|V'_i|\geq \Delta-\numparts>0.99np$, there are still at least $0.04np$ vertices
left in each $V'_i$, $i \geq 3$, that are non-adjacent to both $x$ and $y$. Thus we can apply
Lemma \ref{lem:indep-trans} as above to complete $\{x, y\}$ into an independent
transversal.

The case when $x$ is the vertex of maximum degree has a similar proof but involves one
more step. As in the previous paragraph, we may assume that $x \in V_1$ and has at least
$0.96np$ neighbors in $V_2$, or else we are done. Let $W_2$ be the set of vertices in
$V_2$ that are not adjacent to $x$.  Since $|V_2| = \Delta+1$, we have $W_2 \neq
\emptyset$.  If there exists some $y \in W_2$ that has $< 0.9np$ neighbors in each of the
other $V_i, i\geq 3$, then we can complete $\{x, y\}$ to a full independent transversal as
above.  Otherwise, by Lemma \ref{lem:degree-sequence} the codegree of every pair of
vertices is at most $1.01np^2<0.61np$ and hence each $y \in W_2$ is associated with a
distinct part in which it has $\geq 0.9np$ neighbors.  Yet $x$ has exactly $|W_2| - 1$
neighbors among the other parts $V_i, i\geq 3$, so there must exist $y \in W_2$ such that
$x$ has no neighbors in the part (without loss of generality it is $V_3$) in which $y$ has
$\geq 0.9np$ neighbors. Since $x$ is the unique vertex of maximum degree and $p \gg
\left(\frac{\log^4 n}{n}\right)^{1/3}$, Lemma \ref{lem:degree-sequence} gives
\begin{displaymath}
  d(y) \ \leq \ \Delta-\frac{\sqrt{np}}{\log n} \ < \ \Delta-
  \left\lceil\frac{1}{p}\right\rceil \ \leq \ \Delta-r.
\end{displaymath}
Therefore $V_3$ contains a subset $W_3$ of at least $r+1$ vertices which are not adjacent
to both $x$ and $y$. Since for every $i \geq 4$ at most one vertex in $W_3$ can have more
than $0.81$ neighbors in $V_i$ (by another codegree argument), the pigeonhole principle
ensures that there is a vertex $z \in W_3$ such that $z$ has at most $0.81np$ neighbors in
each $V_i, i\geq 4$.  Also note that $x$ has less than $0.05np$ neighbors in each such
$V_i$, and $y$ has less than $0.11np$. Therefore every $V_i$, $i \geq 4$, has in total
less than $0.05np + 0.11np + 0.81np < (\Delta+1) - 0.03np$ neighbors of any of $\{x, y,
z\}$, so we can apply Lemma \ref{lem:indep-trans} as before to complete $\{x, y, z\}$ into
an independent transversal. \hfill $\Box$

\vspace{0.35cm}

\noindent {\bf Proof of Theorem \ref{thm:main} (ii).}\, We may assume that $p < n^{-1/4}$
because the case $p  \geq n^{-1/4}$ is already a consequence of part (i) of this
theorem.  Fix an arbitrary $\epsilon >0$.  Suppose that $G$ is a graph obtained from
$\gnp$ by adding $(1+\epsilon)\Delta\big\lceil \frac{n}{(1+\epsilon)\Delta}\big\rceil-n$
isolated vertices and $V(G)$ is partitioned into $V_1 \cup \ldots \cup V_r$ with every
$|V_i| = (1+\epsilon)\Delta$. Since $\Delta \geq np$ a.s., we have that $r \leq
\numparts$. We use the same Steps 1--4 to produce a partition of $\cup V_i$ into a
disjoint union of independent transversals. Actually Steps 1--2 can now be made into a
single step, since there is no need here to treat the vertex of maximum degree separately.
The codegree argument implies again that we perform Steps 1--2 at most $r+1$ times.
Moreover, the existence of the independent transversals claimed in these two steps follows
easily from Lemma \ref{lem:indep-trans}. Indeed, suppose that we have deleted
$O\big(\numparts\big)$ independent transversals from $G$.  Since $p \gg \left(\frac{\log
    n}{n}\right)^{1/2}$, we have $1/p=o(np)$ and thus every part still has size at least
$(1+\epsilon/2)\Delta$. Let $x$ be an arbitrary remaining vertex. Since the degree of $x$
is at most $\Delta$, every part still contains at least $\epsilon \Delta/2$ vertices
non-adjacent to $x$.  By Lemma \ref{lem:indep-trans}, we can find an independent
transversal through these vertices which will extend $\{x\}$.

There is no change in the analysis of Step 3 and the same argument as in the proof of part
(i) shows that the total number of transversals deleted from $G$ in Steps 1--3 is at most
$O\big(\numparts \log n\big)$. Since $p \gg \left(\frac{\log n}{n}\right)^{1/2}$, this
number is $o(np)$, and therefore in the beginning of Step 4 each part $V_i$ still has size
$s \geq (1+\epsilon/2)\Delta$.  Recall that in Step 4 we build the remaining $s$ disjoint
independent transversals simultaneously, extending them one vertex at time to cover each
new part $V_{\ell+1}$.  So again we define an auxiliary bipartite graph $H$ whose left
part corresponds to the partial independent transversals $\{T_i\}$ on $V_1, \ldots,
V_\ell$, right part is $V_{\ell+1}$, and the $i$-th vertex on the left is adjacent to $v
\in V_{\ell+1}$ iff $v$ has no neighbors in transversal $T_i$. A perfect matching in $H$
gives a simultaneous extension of each $T_i$.
 
Hence it is enough to verify the Hall condition for $H$, i.e., we must show that for all
$t \leq s$, every set of $t$ vertices on the left has at least $t$ neighbors on the right.
The proof that this holds for all $t \leq s- 40 \numparts \log n$ is exactly the same as
in part (i) and we omit it here.  So suppose that $t > s - 40 \numparts \log n \geq
s-o(np)> (1+\epsilon/3)\Delta$.  Since the degree of every vertex $v \in V_{\ell+1}$ is at
most $\Delta$, it can have neighbors in at most $\Delta<t$ transversals. Therefore there
is at least one transversal in our set of size $t$ which has no neighbors of $v$, and
hence every set of $t> s - 40 \numparts \log n$ vertices on the left has
neighborhood equal to entire right side of $H$.  This verifies the Hall condition and
completes the proof.  \hfill $\Box$

\section{Independent transversals}
\label{sec:third}

In this section, we prove our second theorem.
We only need to consider here the range $\frac{\log^4 n}{n} \ll p \ll \frac{\log^{3/4} n}{\sqrt{n}}$,
since part (ii) of Theorem \ref{thm:main} implies Theorem \ref{thm:main2} for larger values of $p$.
Again, we begin by showing that $\gnp$ satisfies certain properties almost surely.

\subsection{Properties of random graphs}

\begin{lemma}
  \label{lem:degree-sequence-lowprob}
  If $\frac{\log n}{n} \ll p \ll \frac{\log^{3/4} n}{\sqrt{n}}$, then a.s.\ $\gnp$ has the
  following properties:
  \begin{enumerate}
  \item No pair of distinct vertices has more than $3 \log^{3/2} n$ common neighbors.
  \item The maximum degree is strictly between $np$ and $1.01np$.
  \end{enumerate}
\end{lemma}

\noindent 
{\bf Proof.}\,  The codegree $X$ of a
fixed pair of vertices is binomially distributed with parameters $n-2$ and $p^2$. Therefore 
\begin{displaymath}
  \pr{X \geq 3 \log^{3/2} n} \ \leq \ {n-2 \choose 3 \log^{3/2} n} (p^2)^{3 \log^{3/2}
  n} \ \leq \ \left( \frac{enp^2}{3 \log^{3/2} n} \right)^{3 \log^{3/2} n} \ \ll \
(e/3)^{3 \log^{3/2} n} \ = \ o(n^{-2}).
\end{displaymath}
Taking a union bound over all $O(n^2)$ pairs of vertices, we see that
the first property holds a.s.  The second property is a special case of Corollary 3.13
in \cite{Bol}.  \hfill $\Box$

\begin{lemma}
  \label{lem:span-few-edges}
  Let $C\geq 20$ and let $G$ be a graph obtained from the random graph $\gnp$ by
  connecting every vertex to at most $8 \log^2 n$ new neighbors.  Then a.s.\ every subset
  $S \subset V(G)$ of size $|S| \leq C p^{-1} \log^2 n$ spans a subgraph with average
  degree less than $6C\log^2 n$, i.e., contains $< 3C |S| \log^2 n$ edges.
\end{lemma}

\noindent {\bf Proof.}\, Since the edges which we add to the random graph can increase the
number of edges inside $S$ by at most $|S|(8 \log^2 n)/2=4|S|\log^2 n$, it suffices to
show that in $\gnp$ a.s.\ every subset $S$ as above spans less than $eC|S|\log^2 n$ edges.
The probability that this is not the case is at most
\begin{eqnarray*}
  \sum_{m=1}^{Cp^{-1} \log^2 n} {n \choose m} {{m \choose
      2} \choose eCm \log^2 n} p^{eCm \log^2 n} &\leq& \sum_{m=1}^{Cp^{-1} \log^2 n}
  n^m \left(\frac{em}{2eC \log^2 n} \cdot p\right)^{eCm \log^2 n} \\
  &\leq& \sum_{m=1}^{Cp^{-1} \log^2 n} n^m 2^{-eCm \log^2 n} \\
  &\leq& \sum_{m=1}^{Cp^{-1} \log^2 n} \big(n 2^{-eC \log^2 n}\big)^m =o(1),
\end{eqnarray*}
so we are done. \hfill $\Box$

\subsection{Proof of Theorem \ref{thm:main2}}

Fix $\ep > 0$, and suppose we have disjoint subsets $V_1$, \ldots, $V_r$ of $\gnp$, with
all $|V_i| = (1+\ep)\Delta$.  By Lemma \ref{lem:degree-sequence-lowprob}, $r < n/\Delta <
1/p$.  If a vertex $v$ has more than $\frac{\Delta}{\log n}$ neighbors in some $V_i$, say
that $v$ is \emph{locally big}\/ with respect to $V_i$. If it has more than $
\frac{\Delta}{2\log n}$, call it \emph{almost locally big}.  For each $i$, let $B_i$ be
the set of $v$ that are almost locally big with respect to $V_i$.  We claim that $|B_i| <
4 \log n$. Indeed, if $|B_i|\geq 4 \log n$, then Lemma \ref{lem:degree-sequence-lowprob}
together with $\Delta \geq \log^4 n$ and the Jordan-Bonferroni inequality would imply that
the union of neighborhoods in $V_i$ of vertices from $B_i$ is at least
\begin{displaymath}
  (4 \log n) \frac{\Delta}{2\log n} - {4
  \log n \choose 2} 3 \log^{3/2} n \ \geq \ \frac{3}{2}\Delta \ > \ |V_i|,
\end{displaymath}
contradiction.  Next, make each $B_i$ a clique by adding all the missing edges. However,
$\Delta$ will still refer to the maximum degree of the original graph.  Since each vertex
is almost locally big with respect to less than $2\log n$ sets $V_i$, this operation
increases the degree of each vertex by less than $2\log n \cdot 4 \log n=8 \log^2 n\ll
\frac{\Delta}{2\log n}$. Thus every vertex that is locally big after the additions was
almost locally big before.  In particular, there is now an edge between every pair of
vertices that are locally big with respect to the same $V_i$, and there are less than
$r(4\log n) < 4p^{-1} \log n$ locally big vertices in total.

Let $I_1 \subset [r]$ be the set of indices $i$ such that $V_i$ contains more than 
$\frac{\ep}{4} \Delta$ locally big vertices, and define the notation $V_S$ to represent 
$\bigcup_{i \in  S} V_i$.  Note that 
\begin{displaymath}
  |V_{I_1}| \ < \ (1+\ep)\Delta \cdot \left(\frac{\ep}{4} \Delta\right)^{-1} \, 
4p^{-1} \log n \ < \ 20\ep^{-1} p^{-1} \log n
\end{displaymath}
(we can assume here and in the rest of the proof that $\ep$ is sufficiently small).  As
long as there exist $i \not\in I_1$ such that there are more than $(240\ep^{-1}\log^2
n)|V_i|$ crossing edges between $V_i$ and $V_{I_1}$, add $i$ to $I_1$.  Note that each
such index which we add to $V_{I_1}$ increases the number of edges in this set by more than
$(240\ep^{-1}\log^2 n)|V_i|$. Therefore if in this process $I_1$ doubles in size we obtain a set
of size at most $40\ep^{-1} p^{-1} \log n$ with average degree more than
$240\ep^{-1}\log^2 n$, which contradicts Lemma \ref{lem:span-few-edges}. Thus at the end
of the process we have $|I_1| \leq 40\ep^{-1} p^{-1} \log n$.

Given $I_1$, for $t \geq 1$ we recursively define $I_{t+1} \subset I_t$ as follows.  By
Lemma \ref{lem:span-few-edges}, $V_{I_t}$ induces less than $(120\ep^{-1}\log^2
n)|V_{I_t}|$ edges.  Thus, there are less than $2\big(\frac{\Delta}{\log \Delta}\big)^{-1}
\cdot (120\ep^{-1}\log^2 n)|V_{I_t}|$ vertices in $V_{I_t}$ with $> \frac{\Delta}{\log
  \Delta}$ neighbors in this set.  To define $I_{t+1}$ we consider the following process.
Start with $I_{t+1}$ to be the set of all $i \in I_t$ for which $V_i$ has more than
$\frac{\ep}{4} \Delta$ vertices that have $> \frac{\Delta}{\log \Delta}$ neighbors in
$V_{I_t}$. As long as there exist $i \in I_t \setminus I_{t+1}$ such that there are more
than $(240\ep^{-1}\log^2 n)|V_i|$ edges between $V_i$ and $V_{I_{t+1}}$, add $i$ to
$I_{t+1}$.  As above, Lemma \ref{lem:span-few-edges} ensures that this process must stop
before $I_{t+1}$ doubles in size.  Therefore in the end we have
\begin{eqnarray*}
  |I_{t+1}|  &\leq&  2 \left(\frac{\ep}{4}\Delta\right)^{-1} \cdot 2 \left(\frac{\Delta}{\log 
      \Delta}\right)^{-1} \cdot (120\ep^{-1}\log^2 n)|V_{I_t}| \\ 
&\leq&  O\left( \frac{\log^2 n\, \log \Delta}{\Delta^2}|V_{I_t}|\right) \ \leq \ 
O\left( \frac{\log^2 n \, \log \Delta}{\Delta}|I_t|\right)\\
 &\ll& \frac{1}{\log n} |I_t|.
\end{eqnarray*}
Clearly, $|I_1| \leq r \leq n$.  Therefore, when $t \geq \frac{2\log n}{\log\log n}$,
$I_t$ will be empty.  Let $\sigma$ be the smallest index such that $I_\sigma = \emptyset$.  We now
recursively build partial independent transversals $T_\sigma, \ldots, T_1$, where $T_t$ is
an independent transversal on $V_{I_t}$.  Let us say that $T_t$ satisfies property
$\mathbf{P}_t$ if for every $i \not \in I_t$, all the vertices in $T_t$ that
are not locally big with respect to $V_i$ have together at most
$300 (\sigma - t)\frac{\Delta}{\log n}$ neighbors in $V_i$.
It is clear that $T_\sigma = \emptyset$ satisfies $\mathbf{P}_\sigma$, so we can apply the
following lemma inductively to construct $T_1$, an independent transversal on $V_{I_1}$
satisfying $\mathbf{P}_1$.

\begin{lemma}
  \label{lem:first-stage}
  Suppose $t > 1$, and $T_t$ is an independent transversal on $V_{I_t}$ which satisfies
  $\mathbf{P}_t$.  Then we can extend $T_t$ to $T_{t-1}$, an independent transversal on
  $V_{I_{t-1}}$ which satisfies $\mathbf{P}_{t-1}$.
\end{lemma}

We postpone the proof of this lemma until Section \ref{sec:last-lemma}.  Suppose that we
have $T_1$ as described above.  Let $J_1$ be the set of all indices $j \not \in I_1$ such
that some $v \in T_1$ is locally big with respect to $V_j$.  Then, as we did with $I_1$,
as long as there exist $\ell \not\in I_1 \cup J_1$ such that more than $(600\ep^{-1}
\log^2 n)|V_\ell|$ edges cross between $V_\ell$ and $V_{J_1}$, add $\ell$ to $J_1$.  Since
$|T_1|=|I_1|$ and each vertex can be locally big with respect to at most $(1+o(1))\log n$
sets $V_i$, we have that initially $|J_1| \leq (1+o(1))|I_1| \log n \leq 50 \ep^{-1}p^{-1}
\log^2 n$.  Therefore as before, Lemma \ref{lem:span-few-edges} ensures that this process
stops before $J_1$ doubles in size, so the final set $J_1$ has size at most $100
\ep^{-1}p^{-1} \log^2 n$.

As before, we construct a sequence of nested index sets $J_1 \supset \cdots \supset J_\tau
= \emptyset$, where for $t \geq 1$, define $J_{t+1}$ in terms of $J_t$ as follows.  Let
$J_{t+1} \subset J_t$ be the set of all $j \in J_t$ for which $V_j$ contains more than
$\frac{\ep}{4} \Delta$ vertices that have $> \frac{\Delta}{\log \Delta}$ neighbors in
$V_{J_t}$.  Next, as long as there exist $j \in J_t \setminus J_{t+1}$ such that more than
$(600\ep^{-1}\log^2 n)|V_j|$ edges cross between $V_j$ and $V_{J_{t+1}}$, add $j$ to
$J_{t+1}$. Lemma \ref{lem:span-few-edges} again ensures that we stop before $J_{t+1}$
doubles in size, and the same computation as we did for $I_{t+1}$ shows that $|J_{t+1}|
\ll \frac{1}{\log n} |J_t|$. Thus when $t \geq \frac{2 \log n}{\log \log n}$, $J_t$ is
empty. Let $\tau$ be the smallest index for which $J_\tau = \emptyset$.

Next, delete all neighbors of $T_1$ in $V_{J_1}$ and all vertices in $V_{J_1}$ that are
locally big with respect to any $V_k$ with $k \not \in I_1$.  Denote the resulting sets
$V'_j$, $j \in J_1$. We claim that each $V'_j$ still has size at least $\frac{\ep}{2}
\Delta$.  Indeed, at most one $v \in T_1$ can be locally big with respect to $V_j$,
because $T_1$ is an independent set and all vertices that are locally big with respect to
the same part were connected by our construction.  Thus deleting neighbors of this $v$ can
decrease the size of $V_j$ by at most $d(v) < \Delta + 8 \log^2 n =(1+o(1)) \Delta$.  As
for the remaining vertices in $T_1$, which are not locally big with respect to $V_j$,
$\mathbf{P}_1$ ensures that together they have at most $O\big(\sigma \frac{\Delta}{\log
  n}\big) = o(\Delta)$ neighbors in $V_j$, since $\sigma \leq \frac{2 \log n}{\log \log
  n}$.  Also, by construction of $I_1$, every part whose index is not in $I_1$ has at most
$\frac{\ep}{4} \Delta$ locally big vertices.  Hence the size of $V'_j$ is at least $|V_j|
-(1+o(1)) \Delta - \frac{\ep}{4} \Delta \geq \frac{\ep}{2} \Delta$, as claimed.

Let us say that a set $U_t$ satisfies property $\mathbf{Q}_t$ if for every $k \not \in
I_1 \cup J_t$, all the vertices in $U_t$ that
are not locally big with respect to $V_k$ have together at most
$300 (\tau - t)\frac{\Delta}{\log n}$ neighbors in $V_k$.
We need the following  analogue of Lemma \ref{lem:first-stage}.
\begin{lemma}
  \label{lem:second-stage}
  Suppose $t > 1$, and $U_t$ is an independent transversal on $V'_{J_t}$ which satisfies
  $\mathbf{Q}_t$.  Then we can extend $U_t$ to $U_{t-1}$, an independent transversal on
  $V'_{J_{t-1}}$ which satisfies $\mathbf{Q}_{t-1}$.
\end{lemma}
We also postpone the proof of this lemma until Section \ref{sec:last-lemma}.
Starting with $U_\tau = \emptyset$, we iterate this lemma until we obtain $U_1$, an
independent transversal on $V'_{J_1}$ which satisfies $\mathbf{Q}_1$. 
Since $\tau \leq \frac{2 \log n}{\log \log n}$, this property implies that each $V_k$
with $k \not \in I_1 \cup J_1$ has $O\big(\tau \frac{\Delta}{\log n}\big) = o(\Delta)$
vertices with neighbors in $U_1$.

Finally, let $K = [r] \setminus (I_1 \cup J_1)$.  Delete all neighbors of $T_1 \cup U_1$
and all locally big vertices from every $V_k$ with $k\in K$, and denote the resulting sets
by $V'_k$. All $V'_k$ will still have size at least $\big(1 + \frac{\ep}{2}\big) \Delta$,
but now no vertex there has more than $\frac{\Delta}{\log n}$ neighbors in any single set $V_k'$.
Thus, the following result from \cite{LS} implies that for sufficiently large $n$, there
is an independent transversal on $V'_K$, which completes $T_1 \cup U_1$ into an
independent transversal through all parts.

\begin{theorem}
  \label{thm:indep-trans}
  (Loh, Sudakov \cite{LS}) For every $\epsilon>0$ there exists $\gamma>0$ such that the
  following holds.  If $G$ is a graph with maximum degree at most $\Delta$ whose vertex
  set is partitioned into $r$ parts $V_1$, \ldots $V_r$ of size $|V_i| \geq
  (1+\ep)\Delta$, and no vertex has more than $\gamma \Delta$ neighbors in any single part
  $V_i$, then $G$ has an independent transversal.
\end{theorem}

\noindent
This completes the proof of Theorem \ref{thm:main2}, modulo two remaining lemmas. \hfill $\Box$

\subsection{Probabilistic tools}

We take a moment to record two results which we will need for the proofs of the remaining
lemmas.  The first is the symmetric version of the Lov\'asz Local Lemma, which is
typically used to show that with positive probability, no ``bad'' events happen.
\begin{theorem}
  (Lov\'asz Local Lemma \cite{AS}) Let $E_1, \ldots, E_n$ be events. Suppose that there
  exist numbers $p$ and $d$ such that all $\pr{E_i} \leq p$, and each $E_j$ is mutually
  independent of all but at most $d$ of the other events.  If $ep(d+1) \leq 1$, then
  $\pr{\bigcap \overline{E_i}} > 0$.
\end{theorem}
The following result is a short consequence of this lemma, and we sketch its proof for
completeness.
\begin{proposition}
  \label{prop:2e} (Alon \cite{Al-linear-arboricity})
  Let $G$ be a multipartite graph with maximum degree $\Delta$, whose parts $V_1, \ldots,
  V_r$ all have size at least $2e\Delta$.  Then $G$ has an independent transversal.
\end{proposition}
\noindent {\bf Proof.}\, Independently and uniformly select one vertex from each $V_i$,
which we may assume is of size exactly $\lceil 2e\Delta \rceil$.  For each edge $f$ of
$G$, let the event $A_f$ be when both endpoints of $f$ are selected.  The dependencies are
bounded by $2 \lceil 2e\Delta \rceil \Delta -2$, and each $\pr{A_f} \leq \lceil 2e\Delta
\rceil^{-2}$, so the Local Lemma implies this statement immediately.  \hfill $\Box$

\subsection{Proofs of remaining lemmas}
\label{sec:last-lemma}

Since the proofs of Lemmas \ref{lem:first-stage} and \ref{lem:second-stage} are very
similar, we only prove Lemma \ref{lem:first-stage}.  We will simply indicate the two
places where the proofs differ.

\vspace{0.35cm}

\noindent {\bf Proof of Lemma \ref{lem:first-stage}.}\, Fix some $t$ as in the statement
of the lemma.  To extend an independent transversal $T_t$ on the set $V_{I_t}$, satisfying
$\mathbf{P}_t$, to one on the larger set $V_{I_{t-1}}$, satisfying $\mathbf{P}_{t-1}$, we
will use the following key properties of our construction. 
\begin{description}
\item[(i)] For every $i \in I_{t-1} \setminus I_t$, the set $V_i$ contains at most
  $\frac{\ep}{4} \Delta$ vertices that have $> \frac{\Delta}{\log \Delta}$ neighbors in
  $V_{I_{t-1}}$.
\item[(ii)] Each set $V_i$ has size $(1+\ep)\Delta$.
\item[(iii)] For every $i \not \in I_{t-1}$, there are at most $(\beta \log^2 n) |V_i|$
  edges between $V_i$ and $V_{I_{t-1}}$, where we define the constant $\beta$ to be
  $240\ep^{-1}$.
\end{description}
In the case of Lemma \ref{lem:second-stage}, property \textbf{(ii)} is that each set
$V_j'$ has size at least $\frac{\ep}{2} \Delta$, and the constant $\beta$ in property
\textbf{(iii)} is $\beta = 600 \ep^{-1}$.

Let $D = I_{t-1} \setminus I_t$.  From every $V_i$ with $i \in D$, delete all vertices
that have $> \frac{\Delta}{\log \Delta}$ neighbors in $V_{I_{t-1}}$, and all neighbors of
vertices in $T_t$. Denote the resulting sets by $V_i^*$.  Note that now all degrees in the
subgraph on $V_D^* = \bigcup_{i \in D} V_i^*$ are at most $\frac{\Delta}{\log \Delta}$.
Furthermore, we claim that every $|V_i^*| \geq \frac{\ep}{6} \Delta$.  To see this, recall
that at most one vertex $v \in T_t$ can be locally big with respect to $V_i$, because
$T_t$ is independent and all vertices that are locally big with respect to the same part
are connected by our construction. Deleting neighbors of such $v$ can decrease the size of
$V_i$ by at most $d(v) < \Delta + 8\log^2 n =(1+o(1))\Delta$.  The rest of the vertices in
$T_t$ are not locally big with respect to $V_i$, so $\mathbf{P}_t$ implies that they have
less than $O\big(\sigma \frac{\Delta}{\log n}\big) = o(\Delta)$ neighbors in $V_i$ since
$\sigma \leq \frac{2 \log n}{\log \log n}$.  Finally, by property \textbf{(i)} above, in
$V_i$ we will delete at most $\frac{\ep}{4} \Delta$ vertices that have $>
\frac{\Delta}{\log \Delta}$ neighbors in $V_{I_{t-1}}$, so property \textbf{(ii)} implies
that $|V_i^*| \geq (1+\ep)\Delta - (1+o(1))\Delta - \frac{\ep}{4} \Delta \geq
\frac{\ep}{6} \Delta$, as claimed.

In the case of Lemma \ref{lem:second-stage}, recall that by construction all $V_j'$ with
$j \in J_1$ contain no locally big vertices with respect to any part (we deleted all of
them).  Thus, the partial transversal $U_t$ contains no locally big vertices with respect to
$V_j'$.  Property $\mathbf{Q}_t$ then implies that the total number of neighbors that
vertices of $U_t$ have in $V_j'$ is only $O\big(\tau \frac{\Delta}{\log n}\big)=o(\Delta)$.
Hence when we reduce $V_j'$ to $V_j^*$ by deleting all neighbors of $U_t$, and all vertices
that have $> \frac{\Delta}{\log \Delta}$ neighbors in $V_{J_{t-1}}$, the total effect of
$U_t$ is $o(\Delta)$, not $(1+o(1))\Delta$ as above.  Combining this with properties
\textbf{(i)} and \textbf{(ii)}, we see that $|V_j^*| \geq |V_j'| - o(\Delta) -
\frac{\ep}{4}\Delta \geq \frac{\ep}{6}\Delta$, so the claim is still true.  This is the
second and final place in which the proofs of the two lemmas differ, and explains why
Lemma \ref{lem:second-stage} holds with part sizes of only $\frac{\ep}{2} \Delta$, while
Lemma \ref{lem:first-stage} requires part sizes of $(1+\ep)\Delta$.

Returning to the proof of Lemma \ref{lem:first-stage}, randomly select a subset $W_i
\subset V_i^*$ for each $i \in D$ by independently choosing each remaining vertex of $V_i^*$
with probability $\frac{\log^3 \Delta}{\Delta}$, and let $W = \bigcup_{i \in D} W_i$.
Define the following families of bad events.  For each $i \in D$, let $A_i$ be the event
that $|W_i| < \frac{\ep}{8} \log^3 \Delta$, and for each $v \in V_D^*$, let $B_v$ be the
event that $v$ has more than $2 \log^2 \Delta$ neighbors in $W$.  Also, for each $j \not
\in I_{t-1}$, let $C_j$ be the event that the collection of vertices in $W$ that are not
locally big with respect to $V_j$ has neighborhood in $V_j$ of size $> 300
\frac{\Delta}{\log n}$.  We use the Lov\'asz Local Lemma to show that with positive
probability, none of these events happen.

Let us begin by bounding the dependencies.  Say that $A_i$ \emph{lives on}\/ $V_i^*$, $B_v$
\emph{lives on}\/ the neighborhood of $v$ in $V_D^*$, and $C_j$ \emph{lives on}\/ the
neighborhood of $V_j$ in $V_D^*$.  Note that each of our events is completely determined by
the outcomes of the vertices in the set that it lives on.  Hence events living on disjoint
sets are independent.  A routine calculation shows that for any given event, at most
$O(\Delta^3)$ other events can live on sets overlapping with its set; the worst case is
that an event of $C$-type can live on a set that overlaps with the sets of $\leq
(1+\ep)\Delta^3$ other $C$-type events.

It remains to show that each of $\pr{A_i}$, $\pr{B_v}$, and $\pr{C_j}$ are $\ll
\Delta^{-3}$.  The size of $W_i$ is distributed binomially with expectation $\geq \frac{\ep}{6}
\log^3 \Delta$, so by a Chernoff bound, $\pr{A_i} < e^{-\Omega(\log^3 \Delta)} \ll
\Delta^{-3}$.  Similarly, for each $v \in V_D^*$ the expected value of the degree of $v$ in 
$W$ is at most
$\frac{\Delta}{\log \Delta} \cdot \frac{\log^3 \Delta}{\Delta}= \log^2 \Delta$
so $\pr{B_v} < e^{-\Omega(\log^2 \Delta)} \ll \Delta^{-3}$.  For
$\pr{C_j}$, we proceed more carefully.
For each $0 \leq k \leq 8$, let $Y_k$ be the set of vertices in $V_D^*$ that have between
$\frac{\Delta}{\Delta^{(k+1)/8} \log n}$ and $\frac{\Delta}{\Delta^{k/8} \log n}$ many
neighbors in $V_j$.  By property \textbf{(iii)}, the number of edges between $V_{I_{t-1}}$
and $V_j$ is at most $(\beta \log^2 n)|V_j| \leq 2\beta \Delta \log^2 n$.  Therefore,
$|Y_k| \leq 2\beta \Delta^{(k+1)/8} \log^3 n$.  However, 
since $\Delta \geq np \geq \log^4 n$, the probability that at least $30
\Delta^{k/8}$ vertices in $Y_k$ are selected to be in $W$ is bounded by
\begin{eqnarray*}
  \mathbb{P} \ \ \leq \ \
  {2\beta \Delta^{(k+1)/8} \log^3 n \choose 30 \Delta^{k/8}} \left(\frac{\log^3 \Delta}{\Delta}
  \right)^{30 \Delta^{k/8}}
  &\leq& \left( \frac{e \cdot 2\beta \Delta^{1/8} \log^3
      n }{30} \cdot \frac{\log^3 \Delta}{\Delta}\right)^{30\Delta^{k/8}} \\
  &\leq& \left(
    \frac{ e \beta}{15} \cdot \frac{\log^3 \Delta}{\Delta^{1/8}}
  \right)^{30\Delta^{k/8}} \ \ \ll \ \
   \Delta^{-3}.
\end{eqnarray*}
Therefore, with probability $1 - o(\Delta^{-3})$, the collection
of vertices in $W$ that are not locally big with respect to $V_j$ has neighborhood in
$V_j$ of size less than $\sum_{k=0}^8 30 \Delta^{k/8} \frac{\Delta}{\Delta^{k/8} \log n}
< 300 \frac{\Delta}{\log n}$,
and hence $\pr{C_j} \ll \Delta^{-3}$.  

By the Lov\'asz Local Lemma, there exist subsets $W_i \subset V_i^*$ for each $i \in D$ such
that none of the $A_i$, $B_v$, or $C_j$ hold.  In particular, every $|W_i|$ is greater than
$2e$ times the maximum degree in the subgraph induced by $W$, so
Proposition \ref{prop:2e} implies that there exists an independent transversal $T'$ there.
Letting $T_{t-1} = T_t \cup T'$, we obtain an independent transversal on $V_{I_{t-1}}$.
Since $T'\subset W$ and no $C_j$ hold, we have that for every $j \not \in I_{t-1}$, the  
vertices in $T_t \cup T'$ which are not locally big with respect to $V_j$ have together at most
$300 (\sigma - t)\frac{\Delta}{\log n}+300\frac{\Delta}{\log n}=300 (\sigma - (t-1))\frac{\Delta}{\log n}$
neighbors in $V_j$, i.e., $T_t \cup T'$ satisfies $\mathbf{P}_{t-1}$.
\hfill $\Box$

\section{Concluding remarks}

A simple modification of our argument yields a slight improvement of Theorem
\ref{thm:main2}, and shows that the theorem is in fact true for all $p \gg
\frac{\log^{3+\alpha}}{n}$, for any fixed $\alpha > 0$. We decided not to prove
that result here in such generality for the sake of clarity of presentation.
Also, it is not very difficult, using our approach, to prove a statement similar to 
Theorem \ref{thm:main2} for the sparse case, when
$p \sim \frac{c}{n}$ for some constant $c$.  However, these extensions are not as interesting as the
main problem that remains open, which is to study the behavior of the strong
chromatic number of random graphs when  $p \leq n^{-1/2}$.  We are certain that the strong
chromatic number of the random graph $G_{n,p}$ is a.s.\ $(1+o(1))\Delta$ for every 
$p \geq \frac{c}{n}$ for some constant $c$. It would also be very interesting to determine all 
the values of the edge probability $p$ for which almost surely $s\chi(G_{n,p})$ is precisely $\Delta+1$.

\vspace{0.4cm}
\noindent
{\bf Acknowledgments.}\, 
The authors would especially like to thank Bruce Reed for interesting remarks and useful
insights at the early stage of this project.  The idea of studying the strong chromatic
number of random graphs originated from a conversation the second author had with Bruce Reed,
during which it was realized that the strong chromatic number of dense random graphs
should be $\Delta+1$.  We would also like to thank Michael Krivelevich for stimulating
discussions.

\end{document}